\documentclass[12pt,reqno]{amsart}

\title[Weak Centrality of unital $C(X)$-algebras]{Weak Centrality of unital $C(X)$-algebras}

\allowdisplaybreaks[1]

\usepackage{relsize}
\usepackage{amsmath,amsfonts,amssymb,amsthm,mathrsfs}
\usepackage{amsaddr}
\usepackage{hyperref}
\usepackage{cleveref}
\usepackage[margin=2.5cm]{geometry}
\usepackage{setspace}
\usepackage{mathtools}
\usepackage{graphicx}

\usepackage{setspace}
\usepackage{xcolor}

\allowdisplaybreaks[1]

\numberwithin{equation}{section}
\newtheorem{theorem}{\bf Theorem}[section]
\newtheorem{lemma}[theorem]{\bf Lemma}
\newtheorem{cor}[theorem]{\bf Corollary}
\newtheorem{remark}[theorem]{\bf Remark}
\newtheorem{prop}[theorem]{\bf Proposition}

\newtheorem{example}[theorem]{\bf Example}

\newcommand{\Z}{\mathcal{Z}}

\newcommand{\C}{\mathbb{C}}

\newcommand{\FC}{\operatorname{FC}}

\author{Bharat Talwar}
\address{Indian Institute of Technology Bhilai, Dist. Durg, Chhatisgarh - 491002, India}
\author{Prahlad Vaidyanathan}
\address{Indian Institute of Science Education and Research (IISER) Bhopal, Bhauri, Bhopal - 462020, Madhya Pradesh, India.}
\author{Stefan Wagner}
\address{Department of Mathematics and Natural Sciences,
	Blekinge Institute of Technology,
	SE-37179 Karlskrona, Sweden}
\address{Department of Mathematics, 
	Linnaeus University,
	SE-35195  V\"axj\"{o}, Sweden}

\email{btalwar@iitbhilai.ac.in}
\email{prahlad@iiserb.ac.in}
\email{stefan.wagner@bth.se}

\makeatletter
\@namedef{subjclassname@2020}{%
	\textup{2020} Mathematics Subject Classification}
\makeatother



\begin{document}

	\begin{abstract}
		This paper establishes a fibrewise characterization of weak centrality for unital $C(X)$-algebras whose defining homomorphism maps $C(X)$ onto the center: such an algebra is weakly central if and only if each of its nonzero fibres has a unique maximal ideal. 
		This yields a corresponding characterization for arbitrary unital $C^*$-algebras through their canonical fibres over the spectra of their centers. 
		The resulting criterion explains Vesterstr{\o}m's AF-algebra counterexample, whose obstruction is also interpreted through its Bratteli diagram. 
		A parallel criterion characterizes centrality by simplicity of the fibres. 
		Applications to full group $C^*$-algebras give an alternative proof for the discrete Heisenberg group and show that the full group $C^*$-algebra of every countable non-abelian
		torsion-free nilpotent group is not weakly central.
	\end{abstract}

	\maketitle

	\section{Introduction}
	
	The following notation will be used throughout. 
	Given a unital $C^*$-algebra $A$, the symbols $\Z(A)$, $\operatorname{Prim}(A)$, and $\operatorname{Max}(A)$ denote its center, its primitive ideal space, and its maximal ideal space, respectively. 
	The Gelfand spectrum of $\Z(A)$ is denoted by $\Omega_A$, and the maximal ideal corresponding to $z \in \Omega_A$ by $M_z:=\ker(z)$. 
	The Gelfand isomorphism identifies $\Z(A)$ with $C(\Omega_A)$.
	
	\emph{Centrality} asks whether primitive ideals can be recovered from their
	intersections with the center. 
	For each $P \in \operatorname{Prim}(A)$, there is a unique $z \in \Omega_A$ such that $P \cap \Z(A) = M_z$.
	This defines a surjective map
	\[
	\eta: \operatorname{Prim}(A) \to \Omega_A,
	\qquad
	\eta(P) = z.
	\]
	The algebra $A$ is \emph{central} if $\eta$ is injective~\cite{2020IMRN}.
	Further background on central $C^*$-algebras may be found in~\cite[Sec.~9]{Kaplansky}.
	
	A natural weakening of centrality was introduced by Misonou and
	Nakamura~\cite{MN,Misonou}, following ideas of Dixmier~\cite{DixmierJ}.
	Replacing primitive ideals by maximal ideals gives the surjective map
	\[
	\zeta: \operatorname{Max}(A) \to \Omega_A,
	\qquad
	\zeta(M) = z
	\quad
	\Longleftrightarrow
	\quad
	M \cap \Z(A) = M_z.
	\]
	The algebra $A$ is \emph{weakly central} if $\zeta$ is injective~\cite[p.~10]{Somerset}.
	Since $\operatorname{Max}(A) \subseteq\operatorname{Prim}(A)$, every central
	$C^*$-algebra is weakly central.

	Weak centrality admits several equivalent characterizations and is closely related to other structural properties; see, for example,
	\cite{Archbold1976,ArchboldGogicRobert,BresarGogic,PaliwalJain,TalwarJain}.
	The characterization most relevant to the present paper is its equivalence
	to the \emph{center quotient property}: a unital $C^*$-algebra $A$ is
	weakly central if and only if
	\begin{equation}
		\label{eq:clp}
		\Z(A/I)=(\Z(A)+I)/I
	\end{equation}
	for each closed ideal $I\subseteq A$; see
	Vesterstr{\o}m~\cite{Vesterstrom} and
	Archbold--Gogi\'c~\cite[Thm.~3.16]{2020IMRN}.
	
	In an earlier version of the present paper, every AF-algebra was incorrectly
	claimed to be weakly central.
	An anonymous referee pointed out that Vesterstr{\o}m's example~\cite[p.~136]{Vesterstrom} provides a counterexample. 
	The revised paper identifies the mechanism underlying this~failure and places it in a broader bundle-theoretic framework.
	
	The main result gives a fibrewise characterization of weak centrality. 
	More precisely, a unital $C(X)$-algebra whose defining ${}^*$-homomorphism $C(X) \to \Z(A)$ is surjective is weakly central if and only if each of its nonzero fibres has a unique maximal ideal; see~Theorem~\ref{thm:weak-centrality-cx-algebra}. 
	Applying this result to the Gelfand isomorphism $C(\Omega_A) \cong \Z(A)$ shows that a unital~$C^*$-algebra is weakly central if and only if each of its canonical fibres
	over $\Omega_A$ has a unique maximal ideal; see~Corollary~\ref{cor:weak-centrality-canonical-fibres}.
	
	For Vesterstr{\o}m's AF-algebra, the fibrewise criterion identifies an exceptional fibre isomorphic to $\C\oplus\C$ as the obstruction to weak centrality; see~Example~\ref{ex:vesterstrom-bundle-bratteli}. 
	The subsequent Bratteli-diagram analysis interprets this obstruction as the collision of
	distinct quotient labels inside a common matrix summand. 
	It also gives a sufficient separation condition for central lifting; see~Proposition~\ref{prop:eventual-separation}. The unitization $K(H)+\C1$
	shows, however, that the presence of an edge from the quotient part to the
	ideal part does not by itself obstruct weak centrality; see
	Example~\ref{ex:unitized-compacts-bratteli}.
	
	A parallel result characterizes centrality: a unital $C(X)$-algebra with surjective defining homomorphism is central if and only if each of its nonzero fibres is simple; see~Corollary~\ref{cor:centrality-fibres}. 
	Thus, for the canonical fibres of a unital $C^*$-algebra, weak centrality
	corresponds to uniqueness of the maximal ideal, whereas centrality corresponds to simplicity.
	
	The fibrewise criterion also yields applications to full group $C^*$-algebras. 
	It gives an alternative proof that the full group $C^*$-algebra of the discrete Heisenberg group is not weakly central; see~Example~\ref{ex:discrete-heisenberg-group}. 
	More generally, if $G$ is a countable amenable non-abelian discrete group whose $\FC$-center equals $Z(G)$ and for which $G/Z(G)$ is not icc, then $C^*(G)$ is not weakly central; see~Theorem~\ref{thm:large-class-non-weakly-central-groups}. 
	In particular, the full group $C^*$-algebra of every countable non-abelian torsion-free
	nilpotent group is not weakly central; see~Corollary~\ref{cor:torsion-free-nilpotent}.

	\section{Weak centrality: a bundle-theoretic characterization}
	\label{sec:bundle-weak-centrality}
	
	Let $A$ be a unital $C^*$-algebra. 
	In this case, $\operatorname{Max}(A)$ is simply the space of maximal ideals of $A$.
	Denote the inverse Gelfand transform by $C(\Omega_A) \to \Z(A)$, $f \mapsto \widehat{f}$.
	Then, for $z\in\Omega_A$,
	\[
	M_z = \{\widehat f:f\in K_z\},
	\qquad
	K_z := \{f \in C(\Omega_A):f(z)=0\}.
	\]
	
	The following lemma characterizes weak centrality in terms of the ideals of $A$ generated by the maximal ideals $M_z$ of its center.
	
	\begin{lemma}
		\label{lem:weak-centrality-central-ideals}
		Let $A$ be a unital $C^*$-algebra. 
		For $z \in \Omega_A$, let $\langle M_z\rangle$ denote the closed ideal of $A$ generated by $M_z$.
		Then $A$ is weakly central if and only if, for each $z\in\Omega_A$,
		the ideal $\langle M_z\rangle$ is contained in a unique maximal ideal of
		$A$.
	\end{lemma}
	\begin{proof}
		Assume first that $A$ is weakly central, and fix $z \in \Omega_A$.
		By surjectivity of $\zeta$, there~exists $M \in \operatorname{Max}(A)$ such that $M \cap \Z(A) = M_z$.
		Hence $\langle M_z\rangle \subseteq M$. 
		If another maximal ideal $M'$ contains $\langle M_z\rangle$, then
		\[
		M \cap \Z(A) =M' \cap \Z(A) = M_z.
		\]
		Weak centrality gives $M = M'$, proving uniqueness.
		
		Conversely, suppose that $A$ is not weakly central. 
		Then there exist distinct $M_1,M_2\in\operatorname{Max}(A)$ such that
		\[
		M_1 \cap \Z(A)
		=
		M_2 \cap \Z(A)
		=
		M_z
		\]
		for some $z \in \Omega_A$. 
		Therefore $\langle M_z\rangle \subseteq M_1 \cap M_2$, so $\langle M_z\rangle$ is contained in at least two distinct maximal ideals of $A$.
	\end{proof}
	
	The preceding characterization has a natural formulation in terms of $C(X)$-algebras. 
	Let $X$ be a compact Hausdorff space, and let $A$ be a unital $C(X)$-algebra with structure map $\Phi:C(X) \to \Z(A)$, a unital $^*$-homomorphism.
	By Gelfand duality, $\Phi$ determines a continuous map $\sigma_A: \Omega_A \to X$ satisfying
	\[
	\Phi(f) = \widehat{f \circ \sigma_A},
	\qquad f\in C(X).
	\]
	If $\Phi(C(X)) = \Z(A)$, then $\sigma_A$ is a homeomorphism onto its range.

	For $x \in X$, set
	\[
	J_x := \{f\in C(X):f(x)=0\},
	\qquad
	I_x := \overline{\Phi(J_x)A}.
	\]
	The fibre of $A$ over $x$ is
	\[
	A(x) := A/I_x.
	\]
	The family $\{A(x):x\in X\}$ forms an upper semicontinuous field over $X$,
	and $\Phi(f)$ acts on $A(x)$ as multiplication by $f(x)$. 
	The nonzero fibres are precisely those for which $x\in\sigma_A(\Omega_A)$. 
	For background on $C(X)$-algebras, see~\cite[App.~C]{DanaWilliams}.
	
	The following theorem characterizes weak centrality in terms of these fibres.
	
	\begin{theorem}
		\label{thm:weak-centrality-cx-algebra}
		Let $X$ be a compact Hausdorff space and let $A$ be a unital
		$C(X)$-algebra whose structure map $\Phi:C(X) \to \Z(A)$ is surjective.
		Then $A$ is weakly central if and only~if, for each $x \in \sigma_A(\Omega_A)$, the fibre $A(x)$ has a unique maximal ideal.
		
		Whenever these equivalent conditions hold, $\Z(A(x)) \cong \C$ for each $x \in \sigma_A(\Omega_A)$
	\end{theorem}
	\begin{proof}
		By Lemma~\ref{lem:weak-centrality-central-ideals}, $A$ is weakly central if and only if, for each $z\in\Omega_A$, the ideal $\langle M_z\rangle$ is contained in a unique maximal ideal of $A$.
		
		Fix $z \in \Omega_A$ and set $x := \sigma_A(z)$. 
		
		The structure map satisfies $\Phi(J_x) = M_z$.
		Indeed, if $f \in J_x$, then $(f \circ \sigma_A)(z) = 0$, which gives $\Phi(f) \in M_z$.
		Conversely, let $\widehat{g} \in M_z$. 
		Since $\Phi$ is surjective, there exists $f \in C(X)$ such that $\widehat{g} = \Phi(f)$.
		Hence $f(x) = g(z) = 0$, and so $f \in J_x$.
		
		It follows that
		\[
		\langle M_z\rangle
		=
		\overline{M_zA}
		=
		\overline{\Phi(J_x)A}
		=
		I_x.
		\]
		
		By ideal correspondence, $\langle M_z\rangle$ is contained in a unique
		maximal ideal of $A$ if and only if
		\[
		A/\langle M_z\rangle=A/I_x=A(x)
		\]
		has a unique maximal ideal. 
		Since $\sigma_A$ is a homeomorphism onto $\sigma_A(\Omega_A)$, the equivalence~follows.
		
		Finally, assume that the equivalent conditions hold. 
		For $z\in\Omega_A$, set $x := \sigma_A(z)$. 
		The center quotient property~\eqref{eq:clp}, applied to $I_x$, gives
		\[
		\Z(A(x))
		=
		\frac{\Z(A)+I_x}{I_x}
		\cong
		\frac{\Z(A)}{\Z(A) \cap I_x}.
		\]
		The ideal $I_x$ is proper, and $M_z \subseteq \Z(A) \cap I_x \subsetneq \Z(A)$.
		The maximality of $M_z$ therefore gives $\Z(A)\cap I_x = M_z$.
		Consequently,
		\[
		\Z(A(x))
		\cong
		\Z(A)/M_z
		\cong
		\C.
		\qedhere
		\]
	\end{proof}
	
	Theorem~\ref{thm:weak-centrality-cx-algebra} specializes to the canonical $C(\Omega_A)$-algebra structure of
	a unital $C^*$-algebra $A$ as follows:
	
	\begin{cor}
		\label{cor:weak-centrality-canonical-fibres}
		A unital $C^*$-algebra $A$ is weakly central if and only if, for each $z \in \Omega_A$, the
		canonical fibre $A(z)$ has a unique maximal ideal.
	\end{cor}
	
	The following example illustrates how a single fibre with more than one
	maximal ideal obstructs weak centrality.
	
	\begin{example}[Vesterstr{\o}m's AF-algebra]
		\label{ex:vesterstrom-bundle-bratteli}
		Let $X := \mathbb N \cup \{\infty\}$ be the one-point compactification of $\mathbb N$, and consider
		\[
		A := \{a \in C(X,M_2):a(\infty) \in D_2\},
		\]
		where $D_2 \subseteq M_2$ is the diagonal subalgebra. 
		Equivalently, $A$ consists of the convergent~sequences in $M_2$ whose limits are diagonal.
		
		For $n \in \mathbb N$, set
		\[
		A_n := \{a \in A:a(m) = a(\infty)\text{ for all }m>n\}.
		\]
		Then $A_n \cong M_2^{\oplus n} \oplus D_2$ is finite-dimensional. 
		Moreover, $A = \overline{\bigcup_{n \in \mathbb N}A_n}$ since each element of $A$ can be approximated by replacing its tail with its value at $\infty$. 
		Thus $A$ is an AF-algebra.
		
		The algebra $A$ is naturally a unital $C(X)$-algebra with structure map
		$\Phi: C(X) \to \Z(A)$ given by $\Phi(g)(x) := g(x)1_{M_2}$ for $x \in X$.
		
		To determine its center, let $a \in \Z(A)$. 
		For $n \in \mathbb N$, $a(n) \in \Z(M_2) = \C1_{M_2}$, so $a(n) = \lambda_n1_{M_2}$ for some $\lambda_n \in \C$.
		Continuity at $\infty$ implies that $\lambda_n \to \lambda$ and $a(\infty) = \lambda1_{M_2}$ for some $\lambda \in \C$. 
		Thus the function $g \in C(X)$ defined by
		\[
		g(n) := \lambda_n,
		\qquad
		g(\infty) := \lambda,
		\]
		satisfies $a = \Phi(g)$. 
		Consequently, $\Z(A) = \Phi(C(X))$.
		Since $\Phi$ is injective, it is an isomorphism from $C(X)$ onto $\Z(A)$.
		The induced map $\sigma_A: \Omega_A \to X$ is therefore a homeomorphism; after identifying $\Omega_A$ with $X$ in this way, $\sigma_A$ is the identity map.
		
		Evaluation at $x\in X$ identifies the fibres as
		\[
		A(x)\cong
		\begin{cases}
			M_2, & x \in \mathbb N,\\
			D_2 \cong \C \oplus\C, & x = \infty.
		\end{cases}
		\]
		Thus the fibres over $\mathbb N$ are simple, whereas the fibre over
		$\infty$ has two maximal ideals. 
		It follows from Corollary~\ref{cor:weak-centrality-canonical-fibres} that $A$ is not weakly
		central.
		
		This obstruction can also be expressed through the center quotient property~\eqref{eq:clp}.
		The fibre ideal at $\infty$ is
		\[
		I_\infty = \{a \in A:a(\infty) = 0\},
		\]
		and $A/I_\infty \cong D_2 \cong \C \oplus \C$.
		Under the identification $A/I_\infty \cong D_2$, the quotient map is
		evaluation at $\infty$. 
		Since $D_2$ is commutative, $\Z(A/I_\infty) = D_2$.
		The image of $\Z(A)$ under this map consists precisely of the scalar
		matrices in $D_2$. 
		Therefore
		\[
		\frac{\Z(A)+I_\infty}{I_\infty}
		=
		\C1_{D_2}.
		\]
		In particular, the central projection $p := (1,0) \in D_2$ has no central lift to $A$. Consequently,
		\[
		\Z(A/I_\infty)
		\neq
		\frac{\Z(A)+I_\infty}{I_\infty}.
		\]
	\end{example}

	\section{The Bratteli picture of the bundle obstruction}
	\label{sec:bratteli-bundle-obstruction}
	
	The obstruction identified in Example~\ref{ex:vesterstrom-bundle-bratteli} can also be seen in the associated Bratteli~diagram. 
	A central element of the quotient may assign distinct scalar values to different quotient vertices, while the connecting maps may force these values into the same matrix block at a later stage.
	
	To describe this mechanism, let $A = \varinjlim(A_n,\varphi_n)$ be a unital AF-algebra, where
	\[
	A_n = \bigoplus_{v\in V_n}M_{d(v)},
	\]
	and each connecting map $\varphi_n:A_n\to A_{n+1}$ is unital.
	For $v \in V_n$ and $w \in V_{n+1}$, let $m_{w,v}^{(n)}$ denote the multiplicity with which the summand corresponding to $v$ maps into the summand corresponding to $w$.
	Thus $m_{w,v}^{(n)}$ is the number of edges from $v$ to $w$ in the associated
	Bratteli diagram.
	
	Let $I$ be a closed ideal of $A$, and choose a Bratteli presentation adapted
	to $I$. For each $n$, there is then a subset $H_n\subseteq V_n$ such that
	\[
	I_n := I\cap A_n
	=
	\bigoplus_{v\in H_n}M_{d(v)}.
	\]
	Set
	\[
	Q_n := V_n \setminus H_n,
	\qquad
	B_n := A_n/I_n
	\cong
	\bigoplus_{q\in Q_n}M_{d(q)}.
	\]
	After identifying $B_n$ with the complementary direct summand, $A_n = I_n \oplus B_n$.
	
	Ordering the vertices in $H_n$ before those in $Q_n$, the multiplicity
	matrix of $\varphi_n$ has the form
	\[
	\begin{pmatrix}
		P_n & R_n\\
		0   & S_n
	\end{pmatrix}.
	\]
	Here $P_n$ describes the map between the ideal summands, $S_n$ the induced map between the quotient summands, and $R_n$ the edges from $Q_n$ to $H_{n+1}$. 
	The lower-left block is zero because $I$ is an ideal.
	
	The block $R_n$ is a potential source of non-centrality. Its mere
	non-vanishing is not sufficient, however; the relevant phenomenon is the
	collision of distinct quotient labels.
	
	\begin{lemma}[Collision of central labels]
		\label{lem:central-labels}
		With the notation above, let
		\[
		a = \bigoplus_{v \in V_n} \lambda_v1_{d(v)} \in \Z(A_n).
		\]
		If $v_1,v_2 \in V_n$ and $w \in V_{n+1}$ satisfy $m_{w,v_1}^{(n)} \neq 0$ and $m_{w,v_2}^{(n)} \neq 0$, then
		\[
		\operatorname{dist}
		\bigl(\varphi_n(a),\Z(A_{n+1})\bigr)
		\geq
		\frac12|\lambda_{v_1}-\lambda_{v_2}|.
		\]
		In particular, $\varphi_n(a)$ is central if and only if, for each
		$w \in V_{n+1}$, the labels $\lambda_v$ are constant on
		\[
		\{v \in V_n:m_{w,v}^{(n)} \neq 0\}.
		\]
	\end{lemma}
	\begin{proof}
		Fix $w \in V_{n+1}$. 
		Up to unitary conjugacy, the $w$-component of $\varphi_n(a)$ is block diagonal, with scalar blocks $\lambda_v1_{d(v)}$ repeated according to the multiplicities $m_{w,v}^{(n)}$. 
		Its distance from the scalar matrices is at least half the diameter of
		\[
		\{\lambda_v:m_{w,v}^{(n)} \neq 0\}.
		\]
		This gives the asserted estimate. The $w$-component is scalar precisely when
		all labels occurring in it agree, which proves the final assertion.
	\end{proof}
	
	The lemma also applies after telescoping. 
	If $m>n$ and two vertices $v_1,v_2 \in V_n$ both admit paths to a vertex $w \in V_m$, then
	\[
	\operatorname{dist}
	\bigl(\varphi_{m,n}(a),\Z(A_m)\bigr)
	\geq
	\frac12|\lambda_{v_1}-\lambda_{v_2}|,
	\]
	where
	\[
	\varphi_{m,n}
	:=
	\varphi_{m-1} \circ \cdots \circ \varphi_n.
	\]
	
	The obstruction in Vesterstr{\o}m's algebra can now be read directly from
	its Bratteli diagram. Recall from
	Example~\ref{ex:vesterstrom-bundle-bratteli} that $A_n \cong M_2^{\oplus n}\oplus D_2$, with connecting maps
	\[
	\varphi_n(a_1,\ldots,a_n,d)
	=
	(a_1,\ldots,a_n,d,d).
	\]
	In the first copy of $d$ on the right-hand side, the element $d \in D_2$ is regarded as a diagonal matrix in the new $M_2$-summand.
	
	For the ideal $I_\infty$,
	\[
	I_\infty \cap A_n = M_2^{\oplus n} \oplus 0,
	\qquad
	A_n/(I_\infty \cap A_n) \cong D_2.
	\]
	Thus $Q_n$ consists of the two one-dimensional vertices corresponding to the two diagonal entries of $D_2$. 
	Both vertices feed into the new $M_2$-summand in $I_\infty\cap A_{n+1}$.
	
	The projection $p = (1,0) \in D_2$ has labels $1$ and $0$ on these quotient vertices. 
	If $a_n \in \Z(A_n)$ is a central lift of $p$ at level $n$, then Lemma~\ref{lem:central-labels} gives
	\[
	\operatorname{dist}
	\bigl(\varphi_n(a_n),\Z(A_{n+1})\bigr)
	\geq\frac12.
	\]
	Consequently, for every central lift $a_{n+1}\in\Z(A_{n+1})$ of $p$,
	\[
	\|\varphi_n(a_n)-a_{n+1}\|\geq\frac12.
	\]
	The estimate shows that central lifts of $p$ at consecutive finite stages cannot be chosen compatibly. 
	It provides a finite-stage manifestation of the failure of central lifting established in
	Example~\ref{ex:vesterstrom-bundle-bratteli}: the two quotient labels collide inside the new $M_2$-summand.
	
	The absence of edges from the quotient part to the ideal part nevertheless
	provides a useful sufficient condition.
	
	\begin{prop}[A sufficient condition for central lifting]
		\label{prop:eventual-separation}
		Let $A = \varinjlim(A_n,\varphi_n)$ be an AF-algebra and let $I$ be a closed ideal of $A$. Suppose that, after telescoping, the multiplicity~matrices have the form
		\[
		\begin{pmatrix}
			P_n & 0\\
			0   & S_n
		\end{pmatrix}
		\]
		with respect to the decomposition $V_n = H_n \sqcup Q_n$. 
		Then the quotient map $\pi : A \to A/I$ maps $\Z(A)$ onto $\Z(A/I)$. 
		Equivalently,
		\[
		\Z(A/I) = \frac{\Z(A)+I}{I}.
		\]
	\end{prop}
	
	\begin{proof}
		The connecting maps preserve the decompositions $A_n = I_n \oplus B_n$.
		Passing to the inductive limit gives $A \cong I \oplus B$, where $B := \varinjlim(B_n,S_n) \cong A/I$.
		Consequently,
		\[
		\Z(A) \cong \Z(I) \oplus \Z(B),
		\]
		and each element of $\Z(A/I)\cong\Z(B)$ has a central lift of the form
		$0\oplus b$.
	\end{proof}
	
	The condition in Proposition~\ref{prop:eventual-separation} is sufficient
	but not necessary. The following example, suggested by the behaviour of the
	unitization of the compact operators, makes this distinction explicit.
	
	\begin{example}[An edge without a collision]
		\label{ex:unitized-compacts-bratteli}
		Let $H$ be a separable infinite-dimensional Hilbert space and consider $A := K(H)+\C1$.
		Choose strictly increasing finite-rank projections $p_n$ converging strongly
		to $1$ and set
		\[
		A_n := p_nK(H)p_n+\C1
		\cong
		M_{r_n} \oplus \C,
		\qquad
		r_n := \operatorname{rank}(p_n).
		\]
		Then $A = \overline{\bigcup_{n \in \mathbb N}A_n}$, and the connecting maps are given by
		\[
		\varphi_n(a,\lambda)
		=
		\left(
		\operatorname{diag}
		\bigl(a,\lambda1_{r_{n+1}-r_n}\bigr),
		\lambda
		\right).
		\]
		
		For the ideal $I=K(H)$, the matrix vertex belongs to $H_n$ and the scalar vertex belongs to~$Q_n$. 
		At each stage there is an edge from the scalar vertex in $Q_n$ to the matrix vertex in $H_{n+1}$. 
		Thus the block $R_n$ is nonzero for each $n$.
		
		There is, however, only one quotient vertex, so distinct quotient labels
		cannot collide.
		Indeed, $A/I \cong \C$, and each $\lambda \in \C$ has the central lift $\lambda1\in A$. Moreover,
		$K(H)$ is the unique maximal ideal of $A$, and hence $A$ is weakly central.
		
		Thus an edge from $Q_n$ to $H_{n+1}$ does not by itself obstruct weak
		centrality.
	\end{example}
	
	The two examples isolate the relevant Bratteli-diagrammatic mechanism.
	Vanishing of the blocks $R_n$ guarantees central lifting, but is not necessary. 
	The obstruction arises when distinct scalar labels that remain separated in the center of the quotient are forced into a common matrix summand upstairs.

	\section{Centrality and simplicity of fibres}
	\label{sec:centrality-fibres}
	
	The fibrewise characterization of weak centrality has a direct analogue for
	centrality. The difference is that uniqueness of the maximal ideal is
	replaced by simplicity of the fibre.
	
	\begin{theorem}
		\label{thm:centrality-central-ideals}
		Let $A$ be a unital $C^*$-algebra. 
		Then $A$ is central if and only if, for each $z \in \Omega_A$, the ideal $\langle M_z\rangle$ is a maximal ideal of
		$A$.
	\end{theorem}
	\begin{proof}
		Assume first that $A$ is central, and fix $z \in \Omega_A$. 
		By surjectivity of $\zeta$, there exists $Q_z \in \operatorname{Max}(A)$ such that $Q_z \cap \Z(A) = M_z$.
		Hence $\langle M_z\rangle \subseteq Q_z$. 
		Let $P$ be any primitive ideal containing $\langle M_z\rangle$. 
		Then $M_z \subseteq P \cap \Z(A)$.
		The intersection $P \cap \Z(A)$ is a maximal ideal of $\Z(A)$, so $P \cap \Z(A) = M_z$.
		Centrality now gives $P = Q_z$. 
		Thus $Q_z$ is the unique primitive ideal containing $\langle M_z\rangle$. 
		Every closed ideal of a $C^*$-algebra is the intersection of the primitive ideals containing it; hence $\langle M_z\rangle = Q_z$.
		In particular, $\langle M_z\rangle$ is maximal.
		
		The converse can alternatively be obtained by adapting the contrapositive argument in Lemma~\ref{lem:weak-centrality-central-ideals}, with primitive ideals in place of maximal ideals. For completeness, a direct proof is given below:
		
		Suppose that $\langle M_z\rangle$ is maximal for each $z \in \Omega_A$. 
		Let $P_1,P_2 \in \operatorname{Prim}(A)$ satisfy
		\[
		P_1 \cap \Z(A) = P_2 \cap \Z(A).
		\]
		Since $A$ is unital, their common intersection is equal to $M_z$ for some
		$z \in \Omega_A$. 
		Hence $\langle M_z\rangle \subseteq P_i$, $i=1,2$.
		The ideals $P_1$ and $P_2$ are proper, so the maximality of $\langle M_z\rangle$ forces $P_1 = P_2 = \langle M_z\rangle$.
		Therefore $\eta: \operatorname{Prim}(A) \to \Omega_A$ is injective, and $A$ is
		central.
	\end{proof}
	
	For a unital $C(X)$-algebra with surjective structure map, note that the proof of Theorem~\ref{thm:weak-centrality-cx-algebra} gives $\langle M_z\rangle = I_{\sigma_A(z)}$, $z \in \Omega_A$.
	Combining this identity with Theorem~\ref{thm:centrality-central-ideals} yields the following fibrewise characterization.
	
	\begin{cor}
		\label{cor:centrality-fibres}
		Let $X$ be a compact Hausdorff space and let $A$ be a unital $C(X)$-algebra whose structure map $\Phi:C(X) \to \Z(A)$ is surjective. 
		Then $A$ is central if and only if, for each $x \in \sigma_A(\Omega_A)$, the fibre $A(x)$ is simple.
	\end{cor}
	
	In particular, a unital $C^*$-algebra is central if and only if each of its
	canonical fibres over $\Omega_A$ is simple. Thus the canonical fibres
	distinguish the two notions:
	\[
	\begin{array}{ccl}
		\text{weak centrality}
		&\Longleftrightarrow&
		\text{each canonical fibre has a unique maximal ideal},\\[2mm]
		\text{centrality}
		&\Longleftrightarrow&
		\text{each canonical fibre is simple}.
	\end{array}
	\]
	
	For a unital liminal $C^*$-algebra $A$, $\operatorname{Prim}(A)=\operatorname{Max}(A)$; see~\cite[Sec.~4.7.15]{Dixmier}. 
	Consequently, every weakly central unital liminal $C^*$-algebra is central. 
	The following example shows that this implication fails without the liminality assumption.
	
	\begin{example}[Weakly central but not central]
		\label{ex:weakly-central-not-central}
		Let $H$ be a separable infinite-dimensional Hilbert space and set $A := C([0,1],B(H))$.
		Then $\Z(A) = C([0,1])1$, and the canonical fibres are $A(t) \cong B(H)$, $t \in [0,1]$.
		
		The algebra $B(H)$ has a unique maximal ideal, namely $K(H)$, but is not
		simple. 
		Corollary~\ref{cor:weak-centrality-canonical-fibres} therefore shows that $A$ is weakly central, whereas Corollary~\ref{cor:centrality-fibres} shows that it is not central.
		
		For example, the fibre ideal at $0$ is $I_0 = \{f \in A:f(0) = 0\}$, and $A/I_0 \cong B(H)$.
		Thus this fibre has a unique maximal ideal but is not simple.
	\end{example}
	
	Centrality of a unital $C^*$-algebra is equivalent to Hausdorffness of its primitive ideal space; see~\cite[Prop.~3]{Delaroche1967} and \cite[Prop.~2.3.3]{ArchboldThesis}. 
	The analogous statement for weak centrality and~$\operatorname{Max}(A)$ is false.
	
	\begin{example}
		Let
		\[
		A := K(H)+p\C+(1-p)\C\subseteq B(H),
		\]
		where $H$ is a separable infinite-dimensional Hilbert space and $p$ is a projection with infinite-dimensional range and kernel. 
		Then $\operatorname{Max}(A)$ is a two-point Hausdorff space, but $A$ is not
		weakly central; see~\cite[Ex.~3.28]{2020IMRN} and
		\cite[Note~1, p.~257]{DixmierJ}.
	\end{example}

	\section{Applications to group $C^*$-algebras}
	\label{sec:group-cstar-applications}
	
	The fibrewise characterization of weak centrality is particularly useful
	for group $C^*$-algebras whose bundle structures over the spectra of their
	centers can be described explicitly. 
	The discrete Heisenberg group provides the basic example.
	
	\begin{example}[The discrete Heisenberg group]
		\label{ex:discrete-heisenberg-group}
		Let
		\[
		H_3(\mathbb{Z})
		:=
		\left\{
		\begin{pmatrix}
			1 & a & c\\
			0 & 1 & b\\
			0 & 0 & 1
		\end{pmatrix}
		:
		a,b,c\in\mathbb{Z}
		\right\}
		\]
		be the discrete Heisenberg group. 
		Its center is isomorphic to $\mathbb{Z}$, and
		\[
		\Z(C^*(H_3(\mathbb{Z})))
		\cong
		C^*(\mathbb{Z})
		\cong
		C(\mathbb{T}).
		\]
		Thus $C^*(H_3(\mathbb{Z}))$ is naturally a unital $C(\mathbb{T})$-algebra.
		
		For $\theta \in \mathbb{R}$, the fibre over $e^{2\pi i\theta} \in \mathbb{T}$ is
		the quantum 2-torus $A_\theta$, generated by unitaries $U$ and $V$ satisfying $UV = e^{2\pi i\theta}VU$.
		At $\theta = 0$, one has $A_0\cong C(\mathbb T^2)$, which has more than one maximal ideal. Hence Corollary~\ref{cor:weak-centrality-canonical-fibres} implies that $C^*(H_3(\mathbb Z))$ is not weakly central.
	\end{example}
	
	The same argument applies to the free two-step nilpotent groups considered by Omland~\cite{Omland}. 
	Their full group $C^*$-algebras form continuous fields over the duals of their centers, with noncommutative tori as fibres.
	The fibre corresponding to the trivial character of the center is commutative and isomorphic to $C(\mathbb{T}^n)$. 
	Consequently, these group $C^*$-algebras are not weakly central.
	
	The preceding examples are instances of a more general obstruction involving
	the~$\FC$-center. 
	For a discrete group $G$, denote its $\FC$-center by
	\[  
	Z(G)^{\FC}
	:=
	\{g \in G:\text{the conjugacy class of $g$ is finite}\}.
	\]
	Recall that $G$ is \emph{icc} if each nonidentity element has an infinite
	conjugacy class.
	
	\begin{theorem}
		\label{thm:large-class-non-weakly-central-groups}
		Let $G$ be a countable, amenable, non-abelian discrete group such that $Z(G) = Z(G)^{\FC}$ and $G/Z(G)$ is not icc.
		Then $C^*(G)$ is not weakly central.
	\end{theorem}
	\begin{proof}
		Amenability of $G$ gives $C^*(G) = C_r^*(G)$.
		Losert's description of the center of the reduced group $C^*$-algebra~\cite[Thm.~1.1]{Losert}, together with the assumption $Z(G)^{\FC} = Z(G)$, therefore implies that
		\[
		\Z(C^*(G))
		=
		C^*(Z(G))
		\cong
		C\bigl(\widehat{Z(G)}\bigr).
		\]
		Under this identification, $C^*(G)$ is a unital $C\bigl(\widehat{Z(G)}\bigr)$-algebra whose fibre over the trivial~character is $C^*(G/Z(G))$; see~\cite{PackerRaeburn}.
		
		The quotient $G/Z(G)$ is amenable and, by assumption, not icc. Hence
		\[
		C^*(G/Z(G)) = C_r^*(G/Z(G)),
		\]
		and $G/Z(G)$ has a nontrivial finite conjugacy class $\mathcal{C}$.
		Let $\lambda$ denote the left regular representation of $G/Z(G)$ and set
		\[
		a_{\mathcal{C}} := \sum_{s \in \mathcal{C}}\lambda_s.
		\]
		Since $\mathcal{C}$ is invariant under conjugation, $\lambda_t a_{\mathcal C} \lambda_t^* =  a_{\mathcal C}$ for each $t \in G/Z(G)$.
		Thus $a_{\mathcal{C}}$ is a nonscalar central element of $C_r^*(G/Z(G)) = C^*(G/Z(G))$, and consequently
		\[
		\Z(C^*(G/Z(G))) \not \cong \C.
		\]
		By surjectivity of the canonical map from the maximal ideal space of a unital $C^*$-algebra onto that of its center, the fibre $C^*(G/Z(G))$ has more than one maximal ideal. 
		Applying~Corollary~\ref{cor:weak-centrality-canonical-fibres} to this fibre shows
		that $C^*(G)$ is not weakly central.
	\end{proof}
	
	Theorem~\ref{thm:large-class-non-weakly-central-groups} applies in particular
	to torsion-free nilpotent groups.
	
	\begin{cor}
		\label{cor:torsion-free-nilpotent}
		Let $G$ be a countable discrete group that is non-abelian, torsion-free, and
		nilpotent.
		Then $C^*(G)$ is not weakly central.
	\end{cor}
	\begin{proof}
		The group $G$ is amenable, and its $\FC$-center coincides with its center; see~\cite{EckhardtRaum,McLain}. 
		Moreover, the quotient $G/Z(G)$ has nontrivial center because $G$ is non-abelian and nilpotent.
		Any nonidentity element of $Z(G/Z(G))$ has a singleton conjugacy class.
		Hence $G/Z(G)$ is not icc, and the result follows from
		Theorem~\ref{thm:large-class-non-weakly-central-groups}.
	\end{proof}
	
	\begin{example}
		\label{ex:unitriangular-group}
		For $n \geq 3$, let $\operatorname{UT}_n(\mathbb{Z})$ denote the group of
		upper unitriangular $n\times n$ matrices with integer entries. These groups
		are countable, non-abelian, torsion-free, and nilpotent. 
		Hence $C^*(\operatorname{UT}_n(\mathbb{Z}))$ is not weakly central by Corollary~\ref{cor:torsion-free-nilpotent}.
	\end{example}
	
	\begin{example}
		\label{ex:free-nilpotent-groups}
		Let $F_n$ be the free group on $n\geq 2$ generators, and let
		\[
		\gamma_1(F_n):=F_n,
		\qquad
		\gamma_{k+1}(F_n):=[F_n,\gamma_k(F_n)]
		\]
		be its lower central series. 
		For $n,k\geq 2$, the quotient
		\[
		N_{n,k} := F_n/\gamma_{k+1}(F_n)
		\]
		is a countable, non-abelian, torsion-free nilpotent group of class $k$,
		called the free nilpotent group of rank $n$ and class $k$; see~\cite{Omland2020}. 
		Applying Corollary~\ref{cor:torsion-free-nilpotent} to $N_{n,k}$ shows that
		$C^*(N_{n,k})$ is not weakly central.
	\end{example}
	
	The preceding criterion does not cover all group $C^*$-algebras to which
	the fibrewise characterization applies. The following crystallographic
	example has an explicit bundle description and is, moreover, virtually
	abelian.
	
	\begin{example}[A crystallographic group]
		\label{ex:crystallographic-group}
		Let $G$ be the symmetry group of the wallpaper pattern considered
		in~\cite[Sec.~VII.4]{Davidson}. The group contains the translation subgroup
		\[
		\{\tau_{m,n}:(m,n) \in \mathbb{Z}^2\}
		\cong
		\mathbb{Z}^2
		\]
		and a glide reflection $\sigma$. 
		More precisely,
		\[
		G
		=
		\{\tau_{m,n},\sigma\tau_{m,n}:(m,n) \in \mathbb{Z}^2\},
		\]
		with relations
		\[
		\sigma\tau_{m,-n}
		=
		\tau_{m,n} \sigma,
		\qquad
		\sigma^2
		=
		\tau_{1,0}.
		\]
		
		Set $X := \mathbb{T} \times [0,1]$, and  for $z \in \mathbb{T}$, let
		\[
		\mathcal R_z
		:=
		\left\{
		\begin{pmatrix}
			\alpha & z\beta\\
			\beta  & \alpha
		\end{pmatrix}
		:
		\alpha,\beta\in\C
		\right\}
		\subseteq M_2.
		\]
		By~\cite[Cor.~VII.4.2]{Davidson},
		\[
		C^*(G)
		=
		\left\{
		f\in C(X,M_2):
		(\forall (z,t)\in\mathbb T\times\{0,1\})
		\,
		f(z,t) \in \mathcal{R}_z
		\right\}.
		\]
		
		For each $(z,t)\in\mathbb T\times(0,1)$, evaluation at $(z,t)$ maps $C^*(G)$
		onto $M_2$. 
		It follows that every central element of $C^*(G)$ is scalar-valued on $\mathbb{T} \times (0,1)$ and hence, by continuity, on all of $X$.
		Conversely, every scalar-valued section is central. 
		Therefore
		\[
		\Z(C^*(G))
		=
		\{h1_{M_2}:h \in C(X)\}.
		\]
		It follows that $C^*(G)$ is naturally a unital $C(X)$-algebra with surjective structure map~$\Phi: C(X) \to \Z(C^*(G))$ given by $\Phi(h)(x) := h(x)1_{M_2}$ for $x \in X$.
		
		At $x_0:=(1,1)\in X$, a standard cutoff argument shows that evaluation
		identifies the fibre with
		\[
		C^*(G)/I_{x_0}
		\cong
		\mathcal R_1
		=
		\left\{
		\begin{pmatrix}
			a&b\\
			b&a
		\end{pmatrix}
		:
		a,b \in \C
		\right\}.
		\]
		Since it has two maximal ideals, Theorem~\ref{thm:weak-centrality-cx-algebra} implies that $C^*(G)$ is not weakly central.
	\end{example}
	
	\begin{remark}
		The translation subgroup of the group $G$ in
		Example~\ref{ex:crystallographic-group} has index two, so $G$ is virtually
		abelian and of type~I; see~\cite[Sec.~13.11.12]{Dixmier}. 
		Nevertheless, $C^*(G)$ is not weakly central.
	\end{remark}
	
	\section*{Acknowledgments}
	Talwar thanks ANRF: PDF/2023/000688  for their financial support through the NPDF scheme.

\end{document}